\documentclass[10pt]{article}
\include{header}
\usepackage{amsfonts,amssymb,amscd,amsmath,enumerate,verbatim,calc}
\usepackage[english]{babel} 

\usepackage{epsfig}
\usepackage{graphics}
\usepackage{wasysym}
\usepackage{upgreek}
\usepackage{enumitem}
\usepackage{color}
\usepackage{pstricks,pstricks-add,pst-math}
\usepackage{lscape}
\usepackage{stmaryrd}
\usepackage{tikz,tkz-tab}

\setbox0=\hbox{$+$}
\newdimen\plusheight
\plusheight=\ht0
\def\+{\;\lower\plusheight\hbox{$+$}\;}

\setbox0=\hbox{$-$}
\newdimen\minusheight
\minusheight=\ht0
\def\-{\;\lower\minusheight\hbox{$-$}\;}

\setbox0=\hbox{$\cdots$}
\newdimen\cdotsheight
\cdotsheight=\plusheight
\def\cds{\lower\cdotsheight\hbox{$\cdots$}}

\setbox0=\hbox{$+$}

\newtheorem{theorem}{Theorem}[section]
\newtheorem{definition}[theorem]{Definition}

\author{Manjil P. Saikia\\Department of Mathematical Sciences,\\Tezpur University, Napaam\\Dist. - Sonitpur, Pin - 784028\\India\\manjil@gonitsora.com}

\title{Cranks in Ramanujan's Lost Notebook}

\begin{document}
\maketitle

\begin{abstract}
We give a brief review of the work of Ramanujan on cranks that is found in the Lost Notebook. Recent work by Bruce Berndt and his coauthors have brought to light many interesting results of Ramanujan on cranks, which we highlight in this article.
\end{abstract}

\section{Introduction and Motivation}

The motivation for cranks comes from the work of Freeman Dyson \cite{dyson1}, in which he gave combinatorial explanations of the following famous congruences given by Ramanujan

\begin{equation}\label{p1}
 p(5n+4)\equiv 0~(\textup{mod}~5)
\end{equation}

\begin{equation}\label{p2}
 p(7n+5)\equiv 0~(\textup{mod}~7)
\end{equation}

\begin{equation}\label{p3}
 p(11n+6)\equiv 0~(\textup{mod}~11)
\end{equation}

where $p(n)$ is the ordinary partition function defined as follows.

\begin{definition}[Partition Function]
 If $n$ is a positive integer, let $p(n)$ denote the number of unrestricted representations of $n$ as a sum of positive integers, where representations
 with different orders of the same summands are not regarded as distinct. We call $p(n)$ the partition function. 
\end{definition}

In order to give combinatorial explanations of the above, Dyson defined the \textit{rank} 
of a partition to be the largest part minus the number of parts. Let $N(m,t,n)$ denote the number of partitions of $n$ with rank congruent to $m$ modulo $t$. Then Dyson conjectured that
$$N(k,5,5n+4)=\frac{p(5n+4)}{5}, ~~0\leq k \leq 4,$$ and $$N(k,7,7n+5)=\frac{p(7n+5)}{7}, ~~0\leq k \leq 6,$$
which yield combinatorial interpretations of \eqref{p1} and \eqref{p2}.

We record here some of the notations that we shall be using throughout our study. For each nonnegative integer $n$, we set
 
 $$(a)_n:= (a;q)_n:= \prod_{k=0}^{n-1}(1-aq^k), ~~(a)_{\infty}:=(a;q)_{\infty}:=\lim_{n \rightarrow \infty}(a;q)_n, \mid q \mid <1.$$

 We also set $$(a_1, \ldots, a_m;q)_n:=(a_1;q)_n\cdots(a_m;q)_n$$ and $$(a_1, \ldots, a_m;q)_{\infty}:=(a_1;q)_{\infty}\cdots(a_m;q)_{\infty}.$$

The conjectures of Dyson were later proven by Atkin and Swinnerton-Dyer. The generating function for $N(m,n)$, the number of partitions of $n$ with rank $m$ is given by

$$\sum_{m=-\infty}^{\infty}\sum_{n=0}^{\infty}N(m,n)a^mq^n=\sum_{n=0}^{\infty}\frac{q^{n^2}}{(aq;q)_n(q/a;q)_n}.$$ Here $\mid q \mid <1, ~~\mid q \mid < \mid a \mid<1/\mid q \mid.$

 However, the corresponding analogue of the rank doesn't hold for \eqref{p3}, and so Dyson conjectured the existence of another statistic which he called
  \textit{crank}. In his doctoral dissertation Frank G. Garvan defined vector partitions which became the forerunner of the true crank. Then on June 6, 1987 at a 
student dormitory at the University of Illinois, George E. Andrews and Garvan \cite{crank} found the true crank. 

\begin{definition}[Crank]
 For a partition $\pi$, let $\lambda(\pi)$ denote the largest part of $\pi$, let $\mu(\pi)$ denote the number of ones in $\pi$, and let 
 $\nu(\pi)$ denote the number of parts of $\pi$ larger than $\mu(\pi)$. The crank $c(\pi)$ is then defined to be
 
 $$
c(\pi)=\left\{\begin{array}{cc}
{\lambda(\pi)} & if~\mu(\pi)=0,\\
{\nu(\pi)-\mu(\pi)} & if~\mu(\pi)>0.
\end{array}
\right.
$$
\end{definition}

 Let $M(m,n)$ denote the number of partitions of $n$ with crank $m$, and let $M(m,t,n)$ denote the number of partitions of $n$ with crank congruent to 
 $m$ modulo $t$.
 
 For $n\leq 1$ we set $M(0,0)=1, M(m,0)=0$, otherwise $M(0,1)=-1, M(1,1)=M(-1,1)=1$ and $M(m,1)=0$ otherwise.
 
 The generating function for $M(m,n)$ is given by 
 
 \begin{equation}\label{gfc}
\sum_{m=-\infty}^{\infty}\sum_{n=0}^{\infty}M(m,n)a^mq^n=\frac{(q;q)_{\infty}}{(aq;q)_{\infty}(q/a;q)_{\infty}}.  
 \end{equation}

 The crank not only leads to combinatorial interpretations of \eqref{p1} and \eqref{p2}, but also of \eqref{p3}. In fact we have the following result.

 \begin{theorem}[Andrews-Gravan \cite{crank}]
  With $M(m,t,n)$ defined as above, 
  $$M(k,5,5n+4)=\frac{p(5n+4)}{5},~0\leq k \leq 4,$$
  $$M(k,7,7n+5)=\frac{p(7n+5)}{7},~0\leq k \leq 6,$$
  $$M(k,11,11n+6)=\frac{p(11n+6)}{11},~0\leq k \leq 10.$$
 \end{theorem}

\section{Ramanujan and Cranks}

Before giving a brief idea of Ramanujan's contributions on cranks, we define a few things.

 Ramanujan's general theta function $f(a,b)$ is defined by  
 $$f(a,b):=\sum_{n=-\infty}^{\infty}a^{n(n+1)/2}b^{n(n-1)/2}, ~~\mid ab \mid <1.$$

 Our aim is to study the following more general function
 
 $$F_a(q)=\frac{(q;q)_\infty}{(aq;q)_\infty (q/a;q)_\infty}.$$

For brevity, we just 
mention some of the most important ones below. But before that we give the following definition.

\begin{definition}[Dissections]
  If $$P(q):=\sum_{n=0}^{\infty}a_nq^n$$ is any power series, then the m-dissection of $P(q)$ is given by $$P(q)=\sum_{j=0}^{m-1}\sum_{n_j=0}^{\infty}a_{n_jm+j}q^{n_jm+j}.$$
 \end{definition}

The above definition can also be reformulated in the following manner.

\begin{definition}[Dissections]
If $P(q)$ denote any power series in $q$, then the m-dissection of $P$ is given by $$P(q)=\sum_{k=0}^{m-1}q^kP_k(q^m).$$
\end{definition}
 
 In his lost notebook \cite{rama} Ramanujan offers, in various guises, m-dissections for $F_a(q)$ for $m=2,3,5,7,11$.  In particular on page 179
 Ramanujan offers 2- and 3- dissections for $F_a(q)$ in the form of congruences. 
 
 We state the results, without proofs.

 \begin{theorem}[2-dissection]
 We have
 \begin{equation}\label{2d}
F_a(\sqrt{q})\equiv \frac{f(-q^3;-q^5)}{(-q^2;q^2)_{\infty}}+\displaystyle\left(a-1+\frac{1}{a}\displaystyle\right)\sqrt{q}\frac{f(-q,-q^7)}{(-q^2;q^2)_{\infty}}~\displaystyle\left(\textup{mod}~a^2+\frac{1}{a^2}\displaystyle\right). 
 \end{equation}
 \end{theorem}
 
  We note that $\lambda_2=a^2+a^{-2}$, which trivially implies that $a^4\equiv -1~(\textup{mod}~\lambda_2)$ and $a^8\equiv 1~(\textup{mod}~\lambda_2)$. Thus, in \eqref{2d} 
 $a$ behaves like a primitive 8th root of unity modulo $\lambda_2$. If we let $a=\textup{exp}(2\pi i/8)$ and replace $q$ by $q^2$ in the definition of a dissection, \eqref{2d} will give the 2-dissection of $F_a(q)$. In fact, by replacing $q$ with $q^2$, the second definiton of dissection immediately gives us this 2-dissection.

\begin{theorem}[3-dissection]
We have
\begin{align}
 F_a(q^{1/3}) & \equiv \frac{f(-q^2,-q^7)f(-q^4,-q^5)}{(q^9;q^9)_{\infty}}+\displaystyle\left(a-1+1/a\displaystyle\right)q^{1/3}\frac{f(-q,-q^8)f(-q^4,-q^5)}{(q^9;q^9)_{\infty}}\nonumber\\
 & \quad +\displaystyle\left(a^2+\frac{1}{a^2}\displaystyle\right)q^{2/3}\frac{f(-q,-q^8)f(-q^2,-q^7)}{(q^9;q^9)_{\infty}}~\displaystyle\left(\textup{mod}~a^3+1+\frac{1}{a^3}\displaystyle\right).\label{3d}  
\end{align}
\end{theorem}

Again we note that, $\lambda_3=a^3+1+\frac{1}{a^3}$, from which it follows that $a^9\equiv -a^6-a^3\equiv 1~(\textup{mod}~\lambda_3)$. So in \eqref{3d}, $a$ behaves 
 like a primitive 9th root of unity. While if we let $a=\textup{exp}(2\pi i/9)$ and replace $q$ by $q^3$ in the definition of a dissection, \eqref{3d} will give the 3-dissection of $F_a(q)$. In fact, by replacing $q$ with $q^3$, the second definiton of dissection immediately gives us this 3-dissection.

In contrast to \eqref{2d} and \eqref{3d}, Ramanujan offererd the 5-dissection in terms of an equality. 

\begin{theorem}[5-dissection]
Following Ramanujan's notation $f(-q)=f(-q,-q^2)=(q;q)_\infty$, we have
 \begin{align}
 F_a(q) &= \frac{f(-q^2,-q^3)}{f^2(-q,-q^4)}f^2(-q^5)-4\cos^2(2n\pi/5)q^{1/5}\frac{f^2(-q^5)}{f(-q,-q^4)}+2\cos(4n\pi/5)q^{2/5}\frac{f^2(-q^5)}{f(-q^2,-q^3)}\nonumber\\
  &\quad -2\cos(2n\pi/5)q^{3/5}\frac{f(-q,-q^4)}{f^2(-q^2,-q^3)}f^2(-q^5).\label{5d}
 \end{align}

\end{theorem}

 A dissection of a power series is a way to partition that series in classes of similar exponents. We observe that \eqref{5d} has no term with $q^{4/5}$, which is a reflection of \eqref{p1}. This is easy to see, since we get the generating function for the general partition function if we set $a=1$ in $F_a(q)$, and then \eqref{5d} gives us exactly the residues of partitions of $n$ modulo $5$. In fact, one can replace \eqref{5d} by a congruence and in turn \eqref{2d} and \eqref{3d} by equalities. This is done in \cite{bcb1}. In general we can always convert a dissection from one involving roots of unity to one involving a congruence and conversely.
  Ramanujan did not specifically give the 7- and 11-dissections of $F_a(q)$ in \cite{rama}.  However, he vaugely gave some of the coefficients occuring in those dissections. Uniform proofs of these dissections and the others already stated earlier are given in \cite{bcb1}. 
  In \cite{rama}, Ramanujan also recorded various other results related to cranks. For example in pages 179 and 180, he gave tables of $n$ for which $\lambda_n$ satisfies certain congruences.  All these claims are studied
  systematically in \cite{lnb, bcb2} and \cite{bcb3}. Ramanujan also recorded the first 21 coefficients in the power series of $F_a(q)$ in page 58.
  There are many other claims in \cite{rama} that are related to cranks.
  
  \section{{Acknowledgements}}
  
  This paper is one part of the ongoing Masters thesis of the author under the supervision of Prof. Nayandeep Deka Baruah. He thanks Prof. Baruah for his
   valuable guidance and encouragement. The author would like to thank the anonymous referee for valuable comments and suggestions which greatly improved the exposition. The author is supported by INSPIRE Scholarship 422/2009 from DST, Govt. of India.

\end{document}